\pdfoutput=1
\RequirePackage{ifpdf}
\ifpdf 
\documentclass[pdftex]{sigma}
\else
\documentclass{sigma}
\fi

\usepackage{bbm}

\numberwithin{equation}{section}

\newcommand{\Aut}{\operatorname{Aut}}
\newcommand{\Hom}{\operatorname{Hom}}
\newcommand{\Id}{\operatorname{Id}}

\newcommand{\Z}{{\mathbb Z}}
\def\C{{\mathbb C}}

\newcommand{\Oct}{{\mathbb O}}

\newcommand{\Pf}{\operatorname{Pf}}
\newcommand{\SL}{{\mathrm{SL}}}

\newcommand{\SO}{{\mathrm{SO}}}
\newcommand{\Sp}{{\mathrm{Sp}}}
\newcommand{\GG}{{\rm G}_2}

\newcommand{\Tr}{\operatorname{Tr}}

\newtheorem{Theorem}{Theorem}[section]
\newtheorem{Proposition}[Theorem]{Proposition}
 { \theoremstyle{definition}
\newtheorem{Definition}[Theorem]{Definition}
\newtheorem{Remark}[Theorem]{Remark} }

\begin{document}
\allowdisplaybreaks

\newcommand{\arXivNumber}{2009.13573}

\renewcommand{\PaperNumber}{079}

\FirstPageHeading

\ShortArticleName{Triality for Homogeneous Polynomials}

\ArticleName{Triality for Homogeneous Polynomials}

\Author{Laura P.~SCHAPOSNIK~$^{\rm a}$ and Sebastian SCHULZ~$^{\rm b}$}

\AuthorNameForHeading{L.P.~Schaposnik and S.~Schulz}

\Address{$^{\rm a)}$~University of Illinois at Chicago, USA}
\EmailD{\href{mailto:schapos@uic.edu }{schapos@uic.edu}}
\URLaddressD{\url{schapos.people.uic.edu}}

\Address{$^{\rm b)}$~University of Texas at Austin, USA}
\EmailD{\href{mailto:s.schulz@math.utexas.edu }{s.schulz@math.utexas.edu}} 
\URLaddressD{\url{https://web.ma.utexas.edu/users/s.schulz/}}

\ArticleDates{Received February 15, 2021, in final form August 18, 2021; Published online August 27, 2021}

\Abstract{Through the triality of ${\rm SO}(8,\mathbb{C})$, we study three interrelated homogeneous basis of the ring of invariant polynomials of Lie algebras, which give the basis of three Hitchin fibrations, and identify the explicit automorphisms that relate them.}


\Classification{14H60; 31A35; 33C80; 53C07}

\section{Introduction}

This paper is dedicated to the study of the effect of triality on the ring of ad-invariant polynomials on the Lie algebra $\mathfrak{so}(8)$, through the perspective of Hitchin systems. Although this result can be deduced through topological methods (e.g., via the formulae for the Pontrjagin and Euler classes of the spin bundles on an orientable, spinnable 8-manifold \cite{lawson2016spin}) we would like to take here a perspective we have not found elsewhere, and which fits naturally within the study of Higgs bundles.

{\it Triality.} An avatar of triality is triality of vector spaces, which is given by a trilinear form
$\rho\colon V_{1}\times V_{2}\times V_{3}\to \mathbb{R}$ that is non-degenerate in the sense that fixing any two non-zero vectors yields a non-zero linear functional in the third entry. Put differently, fixing a non-zero vector yields a duality of the two remaining vector spaces, i.e., a non-degenerate bilinear form in the usual sense. Vector spaces that are connected via triality can be (non-canonically) identified with a fixed vector space $V$ which is a division algebra. To see this, consider two non-zero vectors $e_1 \in V_1$, and $e_2 \in V_2$. Then, $\rho$ induces isomorphisms $V_2 \xrightarrow{\sim} V_3^\ast$ and $V_1 \xrightarrow{\sim} V_3^\ast$, and thus one can identify these spaces with a vector space $V$. The trilinear form can then be dualized to a map $V \times V \to V$ that we shall call \textit{multiplication}, and the non-degeneracy states precisely that each multiplication has both a left- and a right-inverse, turning $V$ into a division algebra.

The upshot of the above perspective is that triality is a very rigorous phenomenon and over the real numbers it can only appear for vector spaces of dimensions 1, 2, 4 and 8. Across these notes, we are interested in that of dimension 8, where the three vector spaces in question are the vector representation $\Delta_0$ and the two irreducible spin representations $\Delta_1$ and $\Delta_2$ of $\mathrm{Spin}(8)$, all of which are 8-dimensional. The spin representations are self-dual, and so the trilinear form connecting these vector spaces can be seen as the homomorphism $\Delta_0 \times \Delta_1 \to \Delta_2$ that is obtained by restricting the action of the Clifford algebra $\mathrm{Cliff}(8)$ to the vector space $\Delta_0 \simeq \mathbb{R}^8 = \mathrm{Cliff}_1(8)$ of degree 1 elements.
In terms of the trilinear form $\rho$, triality of $\mathrm{Spin}(8)$ means that for every $g \in \mathrm{Spin}(8)$ there exist unique $g_1, g_2 \in \mathrm{Spin}(8)$ such that for all $v_i \in \Delta_i$ one has that
\begin{eqnarray} \rho(v_0, v_1, v_2) = \rho (g v_0, g_1 v_1, g_2 v_2). \label{ecu11}
\end{eqnarray}

{\it Higgs bundles and the Hitchin fibration.} When considering triality between vector spaces and groups, it is natural to ask about its consequences on different mathematical objects defined through those groups and vector spaces. In this paper, we shall ask this question in relation to {\it Higgs bundles}, which were introduced by Hitchin in 1987 for the general linear group \cite{N1}, and whose ``classical'' definition is the following:

 \begin{Definition}\label{def:clasical}\label{clas}
 A {\em Higgs bundle} on a compact Riemann surface $\Sigma$ is a pair $(E,\Phi)$ for a holomorphic vector bundle $E$ on $\Sigma$, and the {\em Higgs field} $\Phi\in H^{0}(\Sigma,\text{End}(E)\otimes K)$, where $K=T^*\Sigma$.
 \end{Definition}

 This definition can be generalized to encompass principal $ G_\mathbb{C}$-bundles, for $ G_\mathbb{C}$ a complex semi-simple Lie group~\cite{N2}, which shall be consider across this paper.
 \begin{Definition}\label{principalLie}\label{defHiggs}\label{complex}
 A {\em $G_\mathbb{C}$-Higgs bundle} on a compact Riemann surface $\Sigma$ is a pair $(P,\Phi)$ where $P$ is a principal $ G_\mathbb{C}$-bundle over $\Sigma$, and the Higgs field $\Phi$ is a holomorphic section of the vector bundle $\operatorname{ad}P\otimes_{\mathbb{C}} K$, for $\operatorname{ad} P$ the vector bundle associated to the adjoint representation and $K=T^*\Sigma$.
 \end{Definition}

When $ G_\mathbb{C}\subset {\rm GL}(n,\mathbb{C})$, a $ G_\mathbb{C}$-Higgs bundle gives rise to a Higgs bundle in the classical sense, with
some extra structure reflecting the definition of $ G_\mathbb{C}$. In particular, classical Higgs bundles are given by ${\rm GL}(n,\mathbb{C})$-Higgs bundles. Through what is known as the non-abelian Hodge correspondence \cite{cor,6,N1,57,yau} and the Riemann--Hilbert correspondence, Higgs bundles manifest themselves as both flat connections and surface group representations, fundamental objects in contemporary mathematics, and closely related to theoretical physics.

By imposing stability conditions, one may form the moduli space $\mathcal{M}_{ G_\mathbb{C}}$ of $G_\mathbb{C}$-Higgs bundles, which in turn has a natural fibration associated to it, the {\it Hitchin fibration}.
 The Hitchin fibration can be defined through a choice of a homogeneous basis $\{p_{i}\}_{i=1}^k$ for the algebra of invariant polynomials of the Lie algebra $\mathfrak{g}_{c}$ of $ G_\mathbb{C}$, of degrees $\{d_{i}\}_{i=1}^k$. Then, the {\em Hitchin fibration}, introduced in \cite{N2}, is given by
\begin{align}
 h\colon \ \mathcal{M}_{ G_\mathbb{C}}&\longrightarrow \mathcal{A}_{ \mathfrak{g}_c}:=\bigoplus_{i=1}^{k}H^{0}\big(\Sigma,K^{d_{i}}\big),\label{base1}\\
 (E,\Phi)&\mapsto (p_{1}(\Phi), \dots, p_{k}(\Phi)).\label{base2}
\end{align}
The map $h$ is referred to as the {\em Hitchin~map}: it is a proper map for any choice of basis of invariant polynomials \cite{N2}, and the space $\mathcal{A}_{ \mathfrak{g}_c}$ is known as the {\it Hitchin base}.\footnote{Notice that the base depends only on the Lie algebra $\mathfrak{g}_c$ as indicated by the notation.}

It is important to note that through the Hitchin fibration, $\mathcal{M}_{ G_\mathbb{C}}$ gives examples of hyperk\"ahler manifolds which are {\it integrable systems} \cite{N2}, leading to remarkable applications in physics.
Moreover, Hausel--Thaddeus \cite{Tamas1} related Higgs bundles to {\it mirror symmetry}, and with Donagi--Pantev presented $\mathcal{M}_{ G_\mathbb{C}}$ as a fundamental example of mirror symmetry for Calabi--Yau manifolds, whose geometry and topology continues to be studied \cite{dopa}. More recently, Kapustin--Witten \cite{Kap} used Higgs bundles and the Hitchin fibration to obtain a physical derivation of the {\it geometric Langlands correspondence} through mirror symmetry. Soon after, Ng\^{o} found the Hitchin fibration a~key ingredients when proving the fundamental lemma in \cite{ngo}.

{\it Summary of our work.} Inspired by the triality induced between three Hitchin fibrations through the triality of Lie groups, Lie algebras and their rings of invariant polynomials, we dedicate this short note to fill a gap we found in the literature when looking for explicit descriptions of correspondences between homogenous bases of the rings of invariant polynomials of Lie algebras arising from the triality of ${\rm SO}(8,\mathbb{C})$. We shall be concerned here with the action of the triality automorphism on the corresponding moduli spaces of Higgs bundles, which has been previously studied by other authors both from a string theory perspective (e.g., see the work of Aganagic--Haouzi--Shakirov \cite{mina} on triality for Coulomb and Higgs branches and related papers) as well as from a mathematics perspective (see the work by Anton Sancho \cite{anton} and Garcia-Prada--Ramanan \cite{oscar}). An automorphism of a Dynkin diagram does not determine a unique lift to an automorphism of the connected, simply-connected complex Lie group it defines. Indeed, it is known that in the case of $\operatorname{Spin}(8, \mathbb{C})$ there are two options up to conjugation by an inner automorphism \cite{wolf-gray}. In particular, these can be chosen so that the fixed locus is either $\mathrm{G}_2$ or $\SL (3, \C )$ (e.g., see~\cite{anton,oscar}).

We will restrict our attention here to a lift corresponding to $\mathrm{G}_2$: using a particular lift $\sigma\colon \mathfrak{so}(8) \to \mathfrak{so}(8)$ of the triality automorphism, we shall study the effect on the base of the Hitchin system explicitly (that is, in a particularly convenient basis). In a different direction, the fixed locus inside the moduli space of Higgs bundles can be described on general grounds via \cite{oscar}. The present work is organized as follows: in Section \ref{sec:triality} we shall give an overview of the group-theoretic construction of triality. Then, in Section \ref{sec:action} we introduce the particular triality automorphism $\sigma$ of $\operatorname{Spin}(8)$, which whilst not difficult to prove, had not been stated in the literature before:

\begin{Proposition}[= Proposition \ref{propf2}] The natural map $\GG \rightarrow \rm Spin(8)$ induced by triality is obtained by combining the action of \begin{equation*}\mathcal{M} := \frac{1}{2} \begin{pmatrix}
-1 & -1 & -1 & -1 \\
\hphantom{-}1 & \hphantom{-}1 & -1 & -1 \\
\hphantom{-}1 & -1 & \hphantom{-}1 & -1 \\
\hphantom{-}1 & -1 & -1 & \hphantom{-}1
\end{pmatrix}
\end{equation*} on all seven quadruples giving the $28$-dimensional Lie algebra $\mathfrak{so}(8)$ in \eqref{quad}, defining an automorphism $\sigma$ of $\mathfrak{so}(8)$ that preserves the Lie bracket, and whose fixed subalgebra is isomorphic to~$\mathfrak{g}_2$, the Lie algebra of $\mathrm{G}_2$.
\end{Proposition}

Our main interest lies in the study of the above action on different homogeneous bases of the ring of invariant polynomials of Lie algebras, since those describe the base of Hitchin fibrations. We hence dedicate Section~\ref{tri:higgs} to study the action of the triality automorphism $\sigma$ described in Proposition~\ref{propf2} on the algebra of invariant polynomials of $\mathfrak{so}(8)$. To this end, recall that a~particular choice of basis is given by the four polynomials $p_1(M) = \Tr\big(M^2\big)$, $p_2 (M) = \Tr\big(M^4\big)$, $p_3 (M) = \Tr\big(M^6\big)$ and $\mathrm{Pf}(M)$, where the latter denotes the Pfaffian. We then prove the following:

\begin{Theorem}[= Theorem \ref{prop:Trafo}] Under the order $3$ automorphism $\sigma$ of $\mathfrak{so}(8)$ in Proposition~{\rm \ref{propf2}}, the basis of $\mathbb{C}[\mathfrak{so}(8)]^{\mathrm{SO}(8)}$ transforms as
\begin{gather*}
\Tr\big(\sigma(M)^2\big) = \Tr\big(M^2\big), \\
\Tr\big(\sigma(M)^4\big) = \frac{3}{8} \Tr\big(M^2\big)^2 - \frac{1}{2} \Tr\big(M^4\big) -12 \Pf (M) , \\
\Pf (\sigma(M)) = - \frac{1}{64} \Tr\big(M^2\big)^2 + \frac{1}{16} \Tr\big(M^4\big) -\frac{1}{2} \Pf (M), \\
\Tr\big(\sigma(M)^6\big) = \frac{15}{64} \Tr\big(M^2\big)^3 - \frac{15}{16} \Tr\big(M^2\big) \cdot \Tr\big(M^4\big) - \frac{15}{2} \Tr\big(M^2\big) \cdot \Pf(M) + \Tr\big(M^6\big).
\end{gather*}
\end{Theorem}

Finally, we conclude the manuscript with some directions of further research for which we envisage the present results shall prove very useful.

\section[Triality of so(8,C)]{Triality of $\boldsymbol{{\mathfrak{so}}(8, \mathbb{C})}$}\label{sec:triality}

 We shall recall here how triality appears for the complex Lie algebra ${\mathfrak{so}}(8, \mathbb{C})$ and the associated simply-connected Lie group $\mathrm{Spin}(8,\mathbb{C})$ from a few different perspectives, which will become useful across these notes. The group $\operatorname{Out}(\mathfrak{g})$ of outer automorphisms of a Lie algebra is the symmetry group of its Dynkin diagram, which for the case at hand is the group $S_3$ of permutations on 3 letters. In particular, these automorphisms permute the three 8-dimensional irreducible representations of $\mathrm{Spin}(8)$ which are given by the vector representation $\Delta_0$ (modelled on $\mathbb{C}^8$) and two chiral spin representations $\Delta_i$ for $i=1,2$.

\subsection{Triality via the octonions}
We shall start by describing the compact real group $\mathrm{G}_2$ as the group of algebra morphisms of the octonions $\Oct$, the maximal real finite-dimensional division algebra. The octonions form a~non-associative, non-commutative unital algebra that is real 8-dimensional. We recall here some of its properties that are needed for our study of triality on Higgs bundles, following \cite{yokota}, which the reader may want to consult for details.

\begin{figure}[h]\centering
\includegraphics[width=0.25\linewidth]{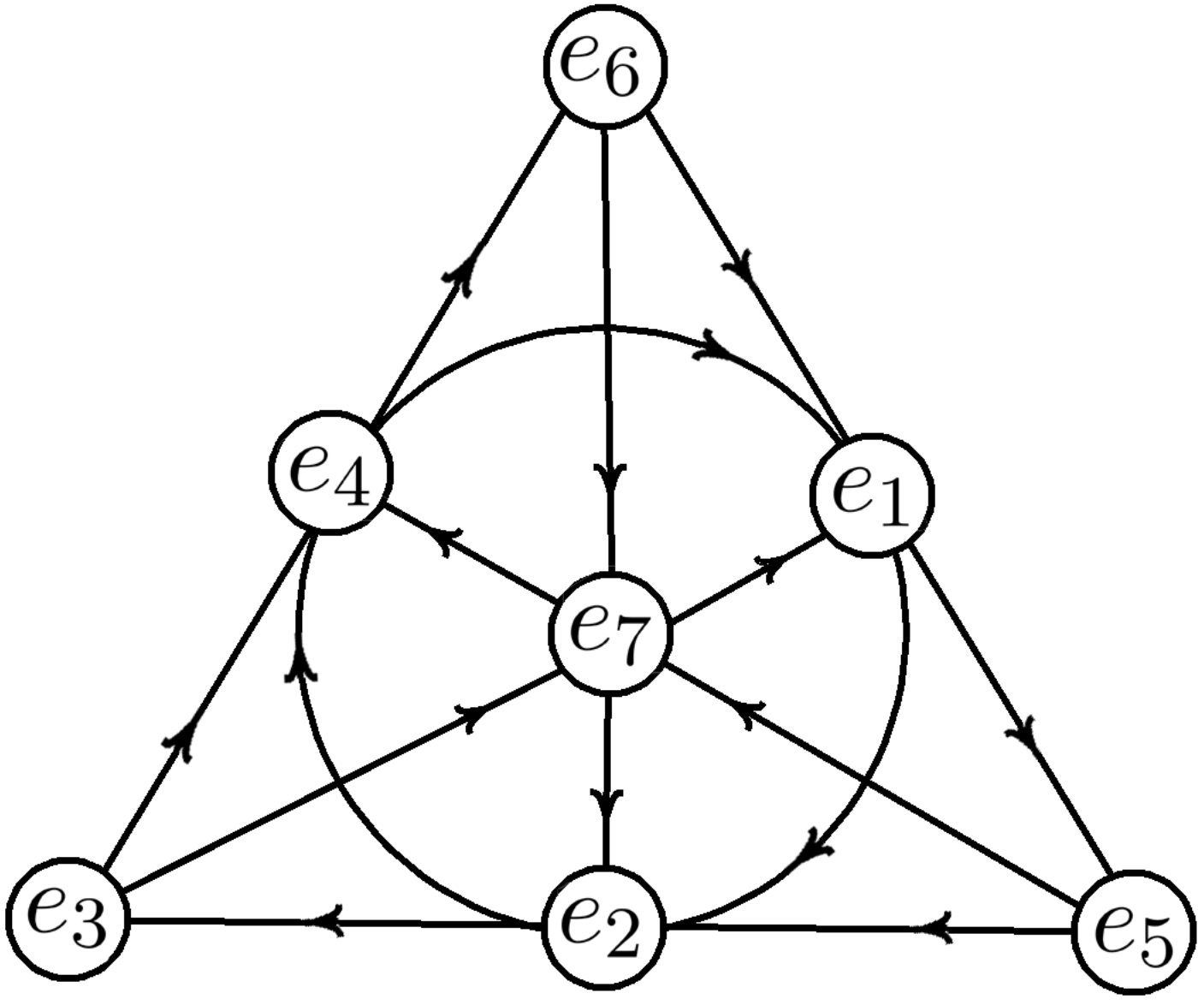}
\caption{The Fano plane captures multiplication of octonions.}\label{pic.fano}
\end{figure}

The starting point is a particular basis $\{e_0=1, e_1,\dots , e_7\}$ for which $e_i^2 = -1$ ($i \neq 0$) and $e_i e_j = -e_j e_i$ ($0 \neq i \neq j \neq 0$). The multiplication of octonions is then completely described by the relations encoded in the Fano plane (see Figure~\ref{pic.fano}). Here, the bottom line for instance reads $e_5 \cdot e_2 = e_3$ and cyclic permutations thereof. Note that $(e_1,e_2,e_4)$ also forms an ordered colinear triple in this way.

The Fano plane is encoding subalgebras: there is a canonical subalgebra isomorphic to $\mathbb{R}$, which is $\mathrm{span}\{e_0\}$. Moreover, every vertex $e_i$ of the diagram identifies a subalgebra $\mathrm{span}\{e_0,e_i\}$ isomorphic to $\mathbb{C}$, and every ordered colinear triple $(e_i, e_j, e_k)$ gives a subalgebra $\mathrm{span}\{e_0,e_i, e_j, e_k\}\!$ isomorphic to the quaternions $\mathbb{H}$. Furthermore, every $e_i$ ($i \neq 0$) sits on exactly three lines, which in the setting of Figure~\ref{pic.fano} are, for indices taken mod 7, given by \begin{equation} \label{eq:QuatEmbedding}
(e_i, e_{i+1}, e_{i+3} ), \qquad (e_i, e_{i+2}, e_{i+6}), \qquad (e_i, e_{i+4}, e_{i+5}).
\end{equation}

Rotating the Fano plane by $2\pi/3$ induces an (order 3) automorphism of $\Oct$ given by $e_0 \mapsto e_0$ and $e_i \mapsto e_{2i}$ where $i \in \{1 , \dots , 7\}$ is taken mod 7. Note that there are also natural order~2 automorphisms given by reflection along one of the central axes, but these have to be accompanied by a sign flip for certain elements to accommodate the correct direction of the arrows.
Similar to the quaternions, for $a_i \in \mathbb{R}$ and $x \in \Oct$, the octonions come equipped with
\begin{itemize}\itemsep=0pt
\item a conjugation $\overline{a_0 + \sum_{i=1}^7 a_i e_i} = a_0 - \sum_{i=1}^7 a_i e_i$,
\item a real part $\operatorname{Re}(x) =\frac{1}{2} (x + \overline{x})$,
\item an inner product $\big(\sum a_\mu e_\mu , \sum b_\nu e_\nu \big) = \sum a_\mu b_\mu = \operatorname{Re}\big((\sum a_\mu e_\mu) \cdot \overline{(\sum b_\nu e_\nu)}\big)$,
\item the induced norm $| x | = \sqrt{(x,x)}$.
\end{itemize}

The group $\mathrm{G}_2$ is the group of algebra automorphisms of the octonions, i.e.,
\begin{equation*}
\mathrm{G}_2 := \{ \alpha \in \Aut_\mathbb{R}(\Oct) \, | \, \alpha (xy) = (\alpha x )(\alpha y) \; \forall x,y \in \Oct \}.
\end{equation*}
The above condition implies in particular that any $\alpha \in \mathrm{G}_2$ obeys $(\alpha x , \alpha y) = (x,y)$ and hence realises $\mathrm{G}_2$ as a closed subgroup of
\begin{equation*}
\mathrm{O}(8) = \mathrm{O}(\Oct) = \{ \alpha \in \Aut_\mathbb{R}(\Oct) \, | \, (\alpha x , \alpha y) = (x,y) \; \forall x,y \in \Oct \}.
\end{equation*}
In particular, $\mathrm{G}_2$ is compact. It is easy to see that it acts trivially on the real part of the octonions (as a morphism of algebras it preserves the unit: $\alpha 1 = 1$), and that the action can be restricted to its orthogonal complement $\Oct' = \textrm{span}\{ e_1 , \dots , e_7 \}$, where $\overline{\alpha e_i} = -\alpha e_i$ ($i=1,\dots,7$), so that $\mathrm{G}_2$ is really a subgroup of
\begin{equation*}
\mathrm{O}(7) = \{ \alpha \in \mathrm{O}(\Oct) \, | \, \alpha 1 = 1 \}.
\end{equation*}

\subsection[Outer automorphisms of Spin(8)]{Outer automorphisms of $\boldsymbol{\mathrm{Spin}(8)}$}
\begin{figure}[h]\centering
\includegraphics[width=0.4\linewidth]{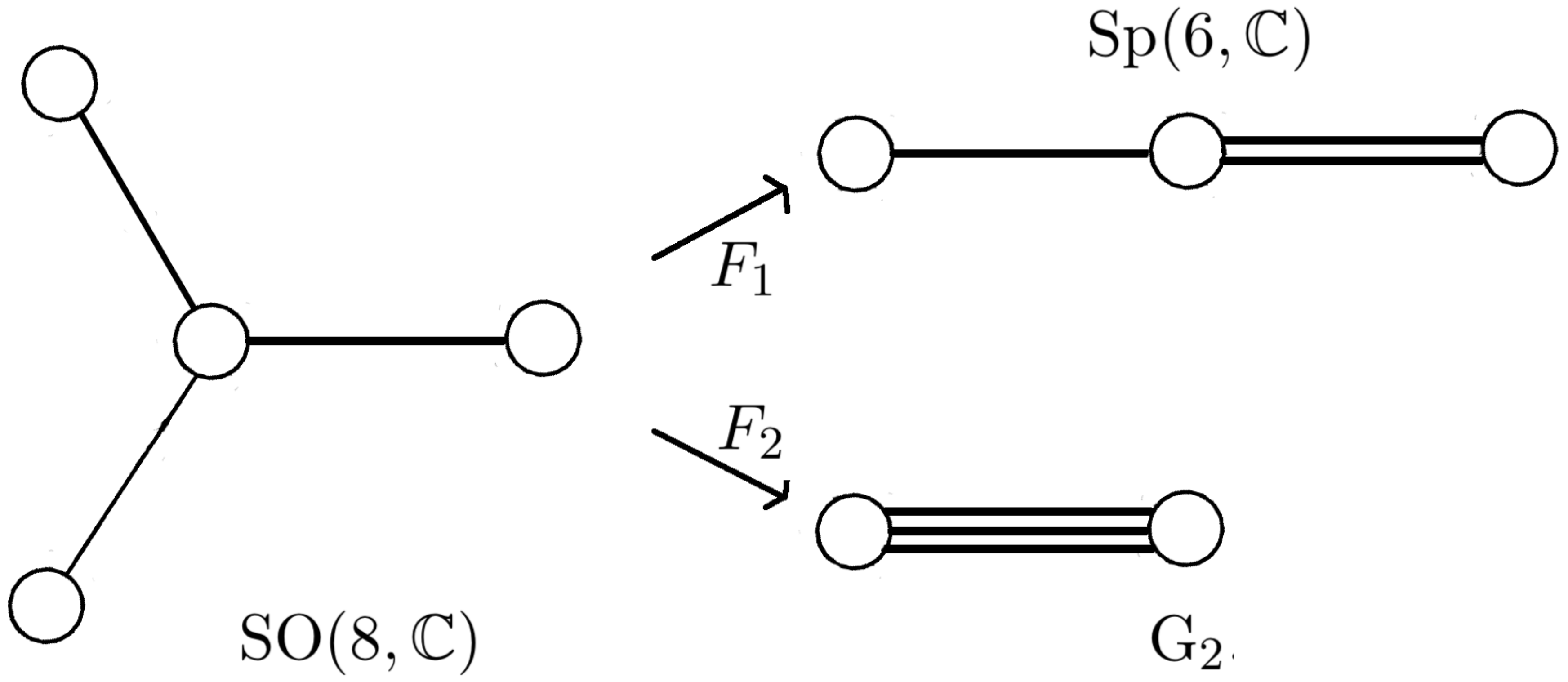}
\caption{The Dynkin diagram $D_4$ presenting its exceptional symmetries, and the two foldings $\mathfrak{f}_1\colon \mathfrak{so}(8) \leadsto \mathfrak{sp}(6)$, and $\mathfrak{f}_2\colon \mathfrak{so}(8) \leadsto \GG$.}\label{folding}
\end{figure}

 The assignment $\sigma_i\colon g \mapsto g_i$ from equation~\eqref{ecu11} is an automorphism that is in fact outer. Recall that an \textit{inner automorphism} of a group $G$ is an automorphism coming from conjugation by some group element $h$, i.e., $g \mapsto h \cdot g \cdot h^{-1} =: C_h(g)$. Inner automorphisms form a normal subgroup $\operatorname{Inn}(G)$ of the group $\operatorname{Aut}(G)$ whose quotient is the group of \textit{outer automorphisms}
 \[\operatorname{Out}(G) := \operatorname{Aut}(G) / \operatorname{Inn}(G).\] Inner automorphisms (by definition) leave the center $Z(G)$ of $G$ invariant and in fact for a simple Lie group $G$, $\operatorname{Inn}(G)$ is naturally isomorphic to $G_{ad} = G/Z(G)$, the adjoint form of the group. If $G$ is additionally simply-connected, $\operatorname{Out}(G)$ is the symmetry group of its associated Dynkin diagram which here is $\operatorname{Out}(\mathrm{Spin}(8)) \simeq S_3$, the group of permutations on 3 elements, as shown in Figure~\ref{folding}.

 Of fundamental importance to us will be the \textit{principle of infinitesimal tri\-a\-li\-ty}~\cite[Theorems~1.3.5 and 1.3.6]{yokota}, which is the infinitesimal version of~(\ref{ecu11}):
\begin{Proposition} \label{prop:InfTriality}
Every $D_1 \in \mathfrak{so} (8)$ determines a unique triple $(D_1, D_2, D_3) \in \mathfrak{so}(8)^3$ such that
\begin{gather*}
 (D_1 x) y + x (D_2 y) = D_3 (xy)
\end{gather*}
for all $x,y \in \Oct$. Furthermore, $D_2 = \sigma (D_1)$, $D_3 = \eta (D_1)$ where $\sigma$ and $\eta$ are outer automorphisms of $\mathfrak{so}(8)$ such that $\sigma^3 = \eta^2 = 1 = (\eta \sigma)^2$. In particular, $\sigma$ and $\eta$ are generators for $S_3 = \operatorname{Out}(\mathfrak{so}(8))$.
\end{Proposition}

The external nodes of the Dynkin diagram correspond to the fundamental representations $\Delta_0$, $\Delta_1$, $\Delta_2$ of $\mathfrak{so}(8)$ (all of which are 8-dimensional), and these are permuted by outer automorphisms, e.g., by $\sigma_1$ and $\sigma_2$ as defined by $\sigma_i\colon g \mapsto g_i$ in equation~\eqref{ecu11}.
 The center of $\mathrm{Spin}(8)$ is $\mathbb{Z}_2\times \Z_2$, which has three elements $\omega_0$, $\omega_1$, $\omega_2$ of order two, such that each $\omega_i$ spans the kernel of~$\Delta_i$. Quotienting ${\rm Spin}(8)$ by one central $\Z_2$ to $\SO(8)$ breaks the $S_3$ symmetry to $\Z_2 \simeq {\rm Out}(\SO (8))$.

 \begin{Remark}
 The above order 2 automorphism of $\SO(8)$ can be represented by conjugation by an element $M \in {\rm O}(8)$ of determinant~$-1$. Conversely, the outer automorphism of ${\rm SO}(8)$ lifts to an automorphism of ${\rm Spin}(8)$ which fixes $\omega_0$ and interchanges $\omega_1$ and $\omega_2$.
 \end{Remark}

\subsection{Folding and fixed point loci}
 The reader should note that the fact that $G_2$ is the fixed point locus of an automorphism of ${\rm Spin}(8)$ that is not inner is no coincidence, as we will explain on the level of Lie algebras in what follows. Recall that any simple, multiply-laced Lie algebra can be realized as the fixed point set of an outer automorphism of a simply-laced Lie algebra. Concretely one has the following:
 \begin{itemize}\itemsep=0pt
 \item The $B$-series $B_n = \mathfrak{so}(2n-1)$ is the fixed point locus of the outer automorphism $X \mapsto \eta^{-1} X \eta$ of $\mathfrak{so}(2n)=D_n$, where $\eta = \operatorname{diag} (-1, 1, 1, \dots, 1) \in {\rm O}(2n)$ has $\det \eta = -1$.
 \item Let $\Omega$ be a non-degenerate, skew-symmetric form of rank $2n$, and recall that $\mathfrak{sp}(2n)$ is isomorphic to the Lie algebra of $(2n \times 2n)$-matrices $X$ for which $\Omega X + X^{\rm T} \Omega = 0$, in other words $X = \Omega^{-1} \big({-}X^{\rm T}\big) \Omega$. The assignment $X \mapsto \Omega^{-1} \big({-}X^{\rm T}\big) \Omega$ defines an outer automorphism of $\mathfrak{sl}(2n) = A_{2n-1}$ that is not inner and whose fixed point locus is the Lie algebra $\mathfrak{sp}(2n) = C_n$.
 \item Finally, $F_4$ is a Lie subalgebra of $E_6$ which is the fixed locus of an outer automorphism, e.g., see~\cite[Section~3.7]{yokota}.
 \end{itemize}

 Of course, different lifts of an outer automorphism to an actual automorphism may have non-isomorphic fixed loci. A more invariant notion is that of a \textit{folding} of a Dynkin diagram of a Lie algebra, which for completion we shall briefly recall here. Outer automorphisms of a Lie algebra $\mathfrak{g}$ are in bijection with symmetries of its underlying Dynkin diagram, hence an outer automorphism $\sigma$ acts on the simple roots $\{ r_i \}_{i \in I}$ via permutation. Let $[i] \in [I]$ denote the set of orbits in $I$ under $\sigma$, and let
 \[ \alpha_{[i]} := \sum_{j \colon [j] = [i]} \alpha_j.\] Then, one can see that the set $\{ \alpha_{[i]} \}_{[i] \in [I]}$ is the set of simple roots for a (typically not simply-laced) Lie algebra $\mathfrak{g}_\sigma$ as long as $\sigma$ does not exchange simple roots that share an edge in the Dynkin diagram (which excludes in particular the outer automorphism of the Dynkin diagram~$A_{2n}$). This new algebra is called the \textit{folded} or \textit{orbit} Lie algebra~\cite{FSS}. In particular, an orbit Lie algebra admits no natural map to the original Lie algebra, but it is Langlands dual to the Lie algebra obtained by taking fixed points. Finally, we should note that it is a standard result that the following foldings occur~\cite{Springer}:
 \begin{itemize}\itemsep=0pt
 \item the series $A_{2n-1}$ folds onto $B_n$ ($n \geq 2$),
 \item the series $D_n$ folds onto $C_{n-1}$ ($n \geq 3$) under an order 2 automorphism,
 \item $E_6$ folds onto $F_4$,
 \item $D_4$ folds onto $G_2$ under the order 3 automorphism.
 \end{itemize}

\section{Triality as an automorphism}\label{sec:action}
In order to understand the appearance of triality via Higgs bundles and the Hitchin fibration, we shall define these subgroups as fixed points of an automorphism to which we turn our attention now, and whose action on the moduli space of Higgs bundles will be studied in the following sections. Recall that the Lie algebra of ${\rm Spin}(8)$ is given by
\begin{equation*}
\mathfrak{so}(8) = \mathfrak{so}(\Oct) = \{ D \in \Hom_\mathbb{R}(\Oct) \, | \, (Dx,y) + (x,Dy) = 0 \; \forall x,y \in \Oct \}
\end{equation*}
with a basis $\{ G_{i,j} \, | \, 0 \leq i < j \leq 7 \}$ defined through
\begin{equation*}
G_{i,j} e_j = e_i, \qquad G_{i,j} e_i = -e_j, \qquad G_{i,j} e_k = 0, \qquad k \neq i,j.
\end{equation*}
It is a 28-dimensional Lie algebra that admits a vector space decomposition into seven 4-dimensional vector spaces with bases
\begin{equation}
\{G_{0,i} , G_{i+1,i+3}, G_{i+2,i+6} , G_{i+4,i+5} \},\label{quad}
\end{equation}
where $i \in \{1, \dots , 7 \}$, and the indices different from $0$ and $7$ are understood mod~7. Notice in particular the resemblance with equation~(\ref{eq:QuatEmbedding}). In this setting, the folding $\mathfrak{f}_1\colon \mathfrak{so}(8) \leadsto \mathfrak{sp}(6)$ is exhibited by taking the fixed locus of an order 2 automorphism, which yields a subalgebra isomorphic to $\mathfrak{so}(7)$. The desired $\mathfrak{sp}(6)$ is then its Langlands dual, and it is in this sense that folding is dual to taking fixed loci.

In order to understand the action of $\mathfrak{f}_2$, consider the linear action of the matrix
\begin{equation}\label{eq:TrafoMatrix}
 \mathcal{M} := \frac{1}{2} \begin{pmatrix}
-1 & -1 & -1 & -1 \\
\hphantom{-}1 & \hphantom{-}1 & -1 & -1 \\
\hphantom{-}1 & -1 & \hphantom{-}1 & -1 \\
\hphantom{-}1 & -1 & -1 & \hphantom{-}1
\end{pmatrix}
\end{equation}
on the four-dimensional subspaces from equation~\eqref{quad}, for which one can show the following:

\begin{Proposition}\label{propf2} The natural inclusion $\GG \hookrightarrow \mathrm{Spin}(8)$ is obtained by combining the action of~$ \mathcal{M}$ from equation~\eqref{eq:TrafoMatrix} on all seven quadruples in equation~\eqref{quad}, which defines an automorphism \begin{equation} \sigma \colon \ \mathfrak{so}(8) \rightarrow \mathfrak{so}(8),\label{autosigma}\end{equation} that preserves the Lie bracket and whose fixed subalgebra is $\mathfrak{so}(8)^\sigma \cong \mathfrak{g}_2$, the Lie algebra of~$\mathrm{G}_2$.
\end{Proposition}
\begin{proof}
It is straightforward to check that $ \mathcal{M}^2 = \mathcal{M}^{\rm T}$ and that $ \mathcal{M}^3 = \mathbbm{1}_4$, i.e., $ \mathcal{M}$ is of order 3 and so is $\sigma$. The fact that it preserves the Lie bracket is a (somewhat tedious, yet straightforward) computation that can be done on the distinguished basis $\{ G_{i,j} \}$ and can be found in~\cite[Lem\-ma~1.3.2]{yokota}.

The proof that the fixed subalgebra is isomorphic to $\mathfrak{g}_2$ is similarly convoluted and can also be found in~\cite[Lemma~1.4.2]{yokota}, where it is again worked out using the particular basis~$\{ G_{i,j} \}$. It crucially uses the \textit{infinitesimal principle of triality}: $D \in \mathfrak{so}(8)$ is fixed by $\sigma$ and $\eta$ from Proposition~\ref{prop:InfTriality} if and only if
\begin{equation*}
 D \in \mathfrak{g}_2 = \textrm{Lie}(G_2) = \mathfrak{der}(\Oct) = \{ D \colon (Dx)y + x(Dy) =D(xy) \; \forall x,y \in \Oct \}.
\end{equation*}
From this perspective, the true difficulty lies in showing that the two definitions of $\sigma$ agree (which can only be done in the basis that the definition requires)~\cite{yokota}.
\end{proof}

\begin{Remark}
It is straightforward to see that $\dim \big( \mathfrak{so}(8)^\sigma \big) =14$, since the $+1$-eigenspace of~$\mathcal{M}$ is 2-dimensional. Using this, together with the facts that $\mathfrak{g}_2 = 14$ and that the fixed locus of a~non-trivial automorphism of $\mathfrak{so}(8)$ is isomorphic to either $\mathfrak{g}_2$ or $\mathfrak{su}(3)$, one might be inclined to use a dimensional argument to conclude the final piece of Proposition~\ref{propf2}. However, one would still need to check that $\sigma$ is indeed an \textit{outer} automorphism, which is in itself a difficult task (even after knowing that $\sigma$ is of order $3$).
\end{Remark}

\begin{Remark}
The principle of infinitesimal triality gives the particular lift $\sigma$ that we consider a very geometric and intrinsic meaning, and it is for this reason that throughout these notes we call $\sigma$ the triality automorphism, even though a different lift from the order-3 automorphism of the underlying Dynkin diagram could rightfully be called ``triality'' as well.
\end{Remark}

\begin{Remark}The folding $\mathfrak{f}_2\colon \mathfrak{so}(8) \leadsto \mathfrak{g}_2$ is obtained by taking the fixed locus of an automorphism, followed by Langlands duality. Moreover, as before, folding does not give rise to a natural map between the two Lie algebras.\end{Remark}

Whilst we have studied above the compact real form $\mathrm{G}_2$, from now on we will care about its complexification which by abuse of notation we will also denote $\mathrm{G}_2$.

\section{Triality and homogeneous invariant polynomials}\label{tri:higgs} Even though the foldings $\mathfrak{f}_1$ and $\mathfrak{f}_2$ do not give rise to natural maps of Lie algebras, they remarkably lead to maps on the level of algebras of invariant polynomials, to which we turn our attention now. For this, we shall first consider how the eigenvalues of matrices are transformed.

\subsection{The choices of homogeneous basis} In order to choose the homogeneous basis of invariant polynomials which we shall be studying, we shall look into how the Hitchin base for different $G_\mathbb{C}$-Hitchin systems are constructed, as described in equation~\eqref{base1}--\eqref{base2}. Since we want to focus on the Lie theoretic aspect of the research here, we shall not go into details on Hitchin systems: the interested reader can find further details on Hitchin base in~\cite{N2} for complex Lie groups, and in~\cite{thesis} for real Lie groups. Moreover, recent applications and open questions in the topic can be found in~\cite{notices}.
 In what follows we take $G_\mathbb{C}$ to be one of the complex Lie groups in Table~\ref{table1} below.

 \begin{table}[h!]\centering\renewcommand{\arraystretch}{1.3}
\begin{tabular}{*{4}{c}}
Lie algebra $\mathfrak{g}$ & Lie group $G_\mathbb{C}$ & Compact real form $\mathfrak{u}$ & dim $\mathfrak{u}$ \\
\hline
$\mathfrak{d}_{4}$ & $\SO(8,\mathbb{C})$	&$\mathfrak{so}(8)$ & $28$ \\
$\mathfrak{b}_{3}$ & $\SO(7,\mathbb{C})$	&$\mathfrak{so}(7)$ & $21$ \\
$\mathfrak{c}_{3}$ & $\Sp(6,\mathbb{C})$	&$\mathfrak{sp}(6)$ & $21$ \\
$\mathfrak{g}_{2}$ & $\GG$	&$\mathfrak{g}_2$ & $14$ \\
 \end{tabular}
\caption{The Lie groups and Lie algebras we consider.}\label{table1}
\end{table}

Since we are looking to further our understanding of the effect of {\it triality} on Higgs bundles, recall that an $\SO(8,\mathbb{C})$-Higgs bundle on a compact Riemann surface~$\Sigma$ is a pair $(E,\Phi)$, for~$E$ a~rk~8 holomorphic vector bundle with a symmetric bilinear form $(\cdot,\cdot)$, and the Higgs field $\Phi\colon E\rightarrow E\otimes T^*\Sigma$, which is a holomorphic map for which $(\Phi v,w)=-(v,\Phi w)$. In local coordinates, $\Phi (z) = M(z) {\rm d}z$ is a holomorphic $\mathfrak{so}(8)$-valued 1-form whose eigenvalues we denote by $\pm \lambda_1$, $\pm \lambda_2$, $\pm \lambda_3$, $\pm \lambda_4$. We shall be interested on how the eigenvalues change under the action induced by the automorphism $\sigma$. For ease of notation we shall denote by $K:=T^*\Sigma$. The characteristic polynomial of the matrix valued map $\Phi$ defines a curve by considering the equation \begin{eqnarray}
\left\{
 \prod_{i=1}^4\big(\eta^2-\lambda_i^2\big)=0
\right\}\subset {\rm Tot}(K).\label{hola1}
\end{eqnarray}
The coefficients $a_i\in H^{0}\big(\Sigma,K^{2i}\big)$ in equation~\eqref{hola1} give a point in the Hitchin base. In order to understand the transformation of this point under triality, it is useful to describe the polynomial in equation~\eqref{hola1} in terms of traces, which we can express as
\begin{eqnarray}
 \eta^8-\left(\sum_{i=1}^4\lambda_i^2\right)\eta^6+\left(\sum_{i<j}\lambda_i^2\lambda_j^2\right)\eta^4 -\left(\sum_{i<j<k}\lambda_i^2\lambda_j^2\lambda_k^2\right)\eta^2+\left(\prod_{i=1}^4\lambda_i^2\right). \label{hola2}
\end{eqnarray}
Since the action in Proposition~\ref{propf2} can be nicely described in terms of actions on traces and Pfaffians, it is useful to describe the characteristic polynomial of equation~\eqref{hola2} in terms of those invariant polynomials, which can be done as follows
\begin{gather*}
\det(\Phi-\eta \cdot \Id)
 = \eta^8-\left(\frac{1}{2}{\rm Tr}\big(\Phi^2\big)\right)\eta^6+
\left(\frac{1}{4}{\rm Tr}\big(\Phi^2\big)^2+\frac{1}{8}{\rm Tr}\big(\Phi^4\big)
\right)\eta^4\nonumber\\
\hphantom{\det(\Phi-\eta \cdot \Id)=}{} +
\left(\frac{1}{48}{\rm Tr}\big(\Phi^2\big)^3-6{\rm Tr}\big(\Phi^2\big){\rm Tr}\big(\Phi^4\big)+8{\rm Tr}\big(\Phi^6\big)
\right)\eta^2
+{\rm Pf}(\Phi)^2. 
\end{gather*}
Hence, a basis of invariant polynomials is given by the Pfaffian $p_4={\rm Pf}(\Phi)$ and
\begin{gather*}
a_1=\frac{1}{2}{\rm Tr}\big(\Phi^2\big), \\
a_2= \frac{1}{4}{\rm Tr}\big(\Phi^2\big)^2+\frac{1}{8}{\rm Tr}\big(\Phi^4\big),
 \\
a_3= \frac{1}{48}{\rm Tr}\big(\Phi^2\big)^3-6{\rm Tr}\big(\Phi^2\big){\rm Tr}\big(\Phi^4\big)+8{\rm Tr}\big(\Phi^6\big).
\end{gather*}

\subsection{The action on the homogeneous basis}
In what follows, we shall consider the values of the invariant polynomials $p_i (M) = \Tr\big(M^{2i}\big)$ for $i=1,2,3$ as well as $p_4(M) = \mathrm{Pf}(M)$ for $\mathfrak{so}(8)$, as well as its Lie subalgebras $\mathfrak{so}(7)$ and $\mathfrak{g}_2$. Recall that any $M \in \mathfrak{so}(8)$ has eigenvalues that come in opposite pairs $\pm \lambda_i$ for $i=1,\dots,4$. For the Lie subalgebras mentioned these restrictions become more severe: It is an easy exercise to see that for $M \in \mathfrak{so}(7)$, written in the 8-dimensional representation obtained by inclusion in $\mathfrak{so}(8)$, two of the eigenvalues must vanish. Moreover, for $M \in \mathfrak{g}_2$ additionally the eigenvalues appear in triples $\lambda_3 = \lambda_1 + \lambda_2$ (for the correct choice of signs). This, together with the subsequent values of the invariant polynomials is given in Tables~\ref{table2}--\ref{table22}, which will be used to describe the induced triality morphism on the homogeneous basis of invariant polynomials for $\mathfrak{so}(8)$.

 \begin{table}[h]\centering\renewcommand{\arraystretch}{1.3}
\begin{tabular}{llccc}
 $\mathfrak{g}$ & Eigenvalues & ${\rm Tr}\big(M^2\big)$& ${\rm Tr}\big(M^4\big)$\\
\hline
$\mathfrak{so}(8)$ & $\pm \lambda_1,\pm \lambda_2,\pm \lambda_3,\pm \lambda_4$	&$2\sum_{i=1}^4\lambda_i^2$ &$2\sum_{i=1}^4\lambda_i^4$ \\
$\mathfrak{so}(7)$ & $\pm \lambda_1,\pm \lambda_2,\pm \lambda_3,0,0$	&$2\sum_{i=1}^3\lambda_i^2$ &$2\sum_{i=1}^3\lambda_i^4$ \\
$\mathfrak{g}_{2}$ & $\pm \lambda_1,\pm \lambda_2,\pm (\lambda_1+\lambda_2),0,0$	& $4\big(\lambda^2_1+\lambda_1\lambda_2+\lambda_2^2\big)$ & $4 \big(\lambda_1 ^2 + \lambda_1 \lambda_2 + \lambda_2^2\big)^2$ \\
\end{tabular}
\caption{Eigenvalues and invariant polynomials, where $\Tr\big(M^4\big) = 1/4 \Tr\big(M^2\big)^2 $ for $M \in \mathfrak{g_2}$.}\label{table2}
\end{table}

 \begin{table}[h]\centering\renewcommand{\arraystretch}{1.3}
\begin{tabular}{llccc}
 $\mathfrak{g}$ & Eigenvalues & ${\rm Tr}\big(M^6\big)$ & ${\rm Pf}(M)$ \\
\hline
$\mathfrak{so}(8)$ & $\pm \lambda_1,\pm \lambda_2,\pm \lambda_3,\pm \lambda_4$ &$2\sum_{i=1}^4\lambda_i^6$	& $\lambda_1\cdot \lambda_2\cdot\lambda_3\cdot\lambda_4$ \\
$\mathfrak{so}(7)$ & $\pm \lambda_1,\pm \lambda_2,\pm \lambda_3,0,0$ &$2\sum_{i=1}^3\lambda_i^6$ 	& 0 \\
$\mathfrak{g}_{2}$ & $\pm \lambda_1,\pm \lambda_2,\pm (\lambda_1+\lambda_2),0,0$ & $ 2\big(\lambda_1^6+\lambda_2^6+(\lambda_1+\lambda_2)^6\big)$ 	& 0
 \end{tabular}
\caption{Eigenvalues and invariant polynomials.} \label{table22}
\end{table}

The algebra of invariant polynomials for $\mathfrak{so}(7)$ is given by $\C [\mathfrak{so}(7)]^{\SO(7)}$ and admits a basis $\{p_1, p_2, p_3 \}$ with $p_i (M) = \Tr\big(M^{2i}\big) $ as before, which gives rise to the natural map
$\C [\mathfrak{so}(7)]^{\SO(7)} \rightarrow \C [\mathfrak{so}(8)]^{\SO(8)}.$
In terms of Hitchin systems, the map $
 \mathcal{A}_{\mathfrak{so}(7)} \rightarrow \mathcal{A}_{\mathfrak{so}(8)}$
is onto the part of the Hitchin base whose preimage under the Hitchin map consists of Higgs bundles with vanishing Pfaffian, or, equivalently, onto the fixed locus under an outer involution (induced by conjugation by a~matrix $A \in \mathrm{O}(8)$ with $\det A = -1$). Recall that a~choice of invariant bilinear form gives an isomorphism $\mathcal{A}_{\mathfrak{sp}(6)} \xrightarrow{\sim} \mathcal{A}_{\mathfrak{so}(7)}$, just like it does for any pair of Langlands dual reductive groups, and the two maps together yield the embedding of the base for the folded Lie algebra.

We shall now turn our attention to the more interesting case of $\mathfrak{g}_2 \hookrightarrow \mathfrak{so}(8)$ as the fixed locus of the triality automorphism. The following theorem establishes how the basis of invariant polynomials transforms:

\begin{Theorem}\label{prop:Trafo}
Under the order $3$ automorphism $\sigma$ of $\mathfrak{so}(8)$ in \eqref{autosigma} induced from equation~\eqref{eq:TrafoMatrix}, the basis $\{p_1 , p_2, p_3, p_4 \}$ of $\mathbb{C}[\mathfrak{so}(8)]^{\mathrm{SO}(8)}$ transforms as
 \begin{gather*}
\Tr\big(\sigma(M)^2\big) = \Tr\big(M^2\big), \\
\Tr\big(\sigma(M)^4\big) = \frac{3}{8} \Tr\big(M^2\big)^2 - \frac{1}{2} \Tr\big(M^4\big) -12 \Pf (M) , \\
\Pf (\sigma(M)) = - \frac{1}{64} \Tr\big(M^2\big)^2 + \frac{1}{16} \Tr\big(M^4\big) -\frac{1}{2} \Pf (M), \\
\Tr\big(\sigma(M)^6\big) = \frac{15}{64} \Tr\big(M^2\big)^3 - \frac{15}{16} \Tr\big(M^2\big) \cdot \Tr\big(M^4\big) - \frac{15}{2} \Tr\big(M^2\big) \cdot \Pf(M) + \Tr\big(M^6\big).
\end{gather*}
\end{Theorem}

\begin{proof}The first identity is straightforward: The space of invariant polynomials of degree~2 is one-dimensional, and $\sigma$ defines an action of $\mathbb{Z}_3$ on it, hence acts through multiplication by a cubic root of unity. It is easy to see from the definition of $\sigma$ that it acts purely real, hence leaving $\Tr\big(M^2\big)$ invariant.

For the other three polynomials one need to perform some further analysis. In order to understand the action of~$\sigma$, we shall consider the values of the invariant polynomials in terms of the matrix entries of $M=\{M_{ij}\}$. As for any antisymmetric matrix, the trace of its square is
\begin{gather*}
\Tr\big(M^2\big) = \sum_i \big(M^2\big)_{ii} = \sum_{i,j} M_{ij} M_{ji} = \sum_{i,j} -(M_{ij})^2 = -2 \sum_{i<j} M_{ij}^2,
\end{gather*}
where $i,j,\ldots = 0,\dots,7$ unless otherwise noted. Since $M^2$ is itself symmetric, then
\begin{gather*}
\Tr\big(M^4\big) = \sum_{i,j} \big(M^2\big)_{ij}^2 = 2 \sum_{i<j} M_{ij}^4 + 4 \sum_{i<j<k} \big( M_{ij}^2 M_{ik}^2 + M_{ij}^2 M_{jk}^2 + M_{ik}^2 M_{jk}^2 \big) \\
\hphantom{\Tr\big(M^4\big) = \sum_{i,j} \big(M^2\big)_{ij}^2 =}{} + 8 \!\!\sum_{i<j<k<l} \!\!\big( M_{ij} M_{ik} M_{jl} M_{kl} - M_{ij} M_{jk} M_{kl} M_{il} + M_{ik} M_{il} M_{jk} M_{jl} \big).
\end{gather*}

The expression for $\Tr\big(M^6\big)$ is increasingly complicated, but can be calculated in a similar fashion.
Lastly, we can compute the Pfaffian from the expression of the determinant, see Figure~\ref{pic:PfaffM}, to obtain
\begin{gather}
 \Pf (M) = \frac{1}{6 \cdot 8}\sum_{\eta \in S_7} \mathrm{sgn}(\eta) \cdot M_{0 \eta(1)} \cdot M_{\eta(2) \eta(3)} \cdot M_{\eta(4) \eta(5)} \cdot M_{\eta(6) \eta(7)} \nonumber\\
\hphantom{\Pf (M)}{}= \frac{1}{4! \cdot 2^4}\sum_{\eta \in S_8} \mathrm{sgn}(\eta) \cdot M_{\eta(0) \eta(1)} \cdot M_{\eta(2) \eta(3)} \cdot M_{\eta(4) \eta(5)} \cdot M_{\eta(6) \eta(7)},
\label{pref}
\end{gather}
where $S_7$ (resp.\ $S_8$) is the symmetric group on the letters $\{1, \dots ,7\}$ (resp.\ on $\{0, \dots ,7 \}$).

\begin{figure}[h!]\centering
\includegraphics[width=0.9\linewidth]{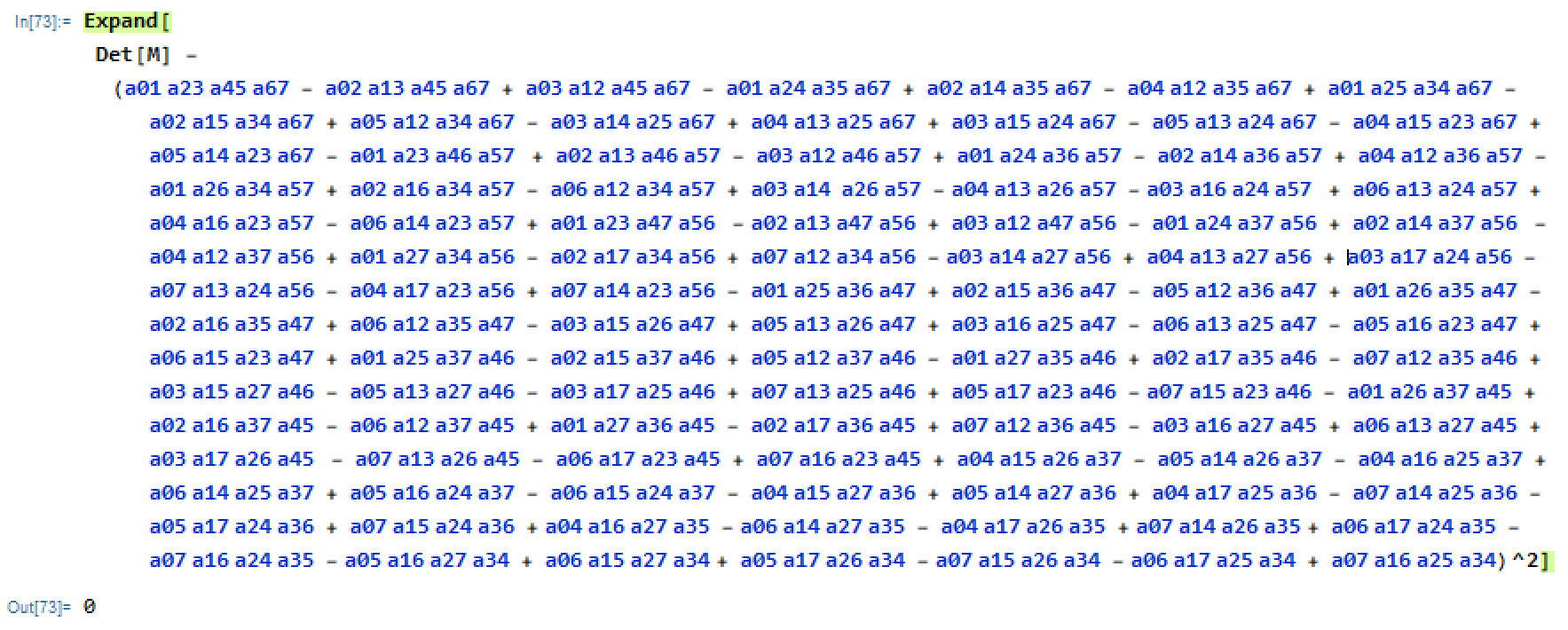}
\caption{The Pfaffian $\mathrm{Pf}(M) = \sqrt{\det M}$.}\label{pic:PfaffM}
\end{figure}

One should note that the prefactor in the first line of equation~\eqref{pref} arises
(compared to Figure~\ref{pic:PfaffM})
 from permuting the individual factors without the subscript $0$ (alternatively, from imposing $\eta(2) < \eta(4) < \eta(6)$), as well as from ordering the individual subscripts by size using $M_{ij} = -M_{ji}$ (alternatively by imposing $\eta(i) < \eta(i+1)$ for $i=2,4,6$), and similarly for the second line.

\looseness=-1 Finally, recall that the automorphism $\sigma$ acting on $\mathfrak{so}(8)$ is induced from the linear action of the matrix defined in~equation~(\ref{eq:TrafoMatrix}) on the linear subspaces spanned for $i \in \{1,\dots,7\} \mod 7$ by
\[
\{M_{0,i}, M_{i+1,i+3}, M_{i+2,i+6}, M_{i+4,i+5} \}.
\]

As an example, the first column of the transformed matrix $X:= \sigma(M)$ takes the form
\[\frac{1}{2}\left(
\begin{matrix}
 0 \\
 M_{01} + M_{24} + M_{37} + M_{56} \\
 M_{02} + M_{35} + M_{67} - M_{14} \\
 M_{03} + M_{46} - M_{17} - M_{25} \\
 M_{04} + M_{12} + M_{57} - M_{36} \\
 M_{05} + M_{23} - M_{16} - M_{47} \\
 M_{06} + M_{15} + M_{34} - M_{27} \\
 M_{07} + M_{13} + M_{26} + M_{45}
\end{matrix}
\right).
\]
By degree reasons, $ \Tr \big(X^4\big)$ can be expressed as \begin{equation*}
 \Tr \big(X^4\big) = A \cdot \Tr \big(M^2\big)^2 + B \cdot \Tr \big(M^4\big) + C \cdot \Pf (M),
\end{equation*}
for some constants $A$, $B$, $C$ which we shall determine next. To this end, note that $\Tr \big(X^4\big)$ has the following shape
\begin{gather*}
 \Tr\big(X^4\big)= \frac{1}{2} \big( M_{01}^4 + M_{02}^4 + \cdots \big)
 +\big( M_{01}^2 M_{02}^2 + M_{01}^2 M_{12}^2 + M_{02}^2 M_{12}^2 + \cdots\big) \\
\hphantom{\Tr\big(X^4\big)=}{}
+ 3 \big( M_{01}^2 M_{23}^2 + M_{01}^2 M_{24}^2 + \cdots \big)
 + 4 \big( M_{01} M_{12} M_{23} M_{03} - M_{02} M_{03} M_{12} M_{13} \pm \cdots \big) \\
\hphantom{\Tr\big(X^4\big)=}{} -12 \big( M_{01} M_{23} M_{45} M_{67} - M_{02} M_{13} M_{45} M_{67} \pm \cdots \big).
\end{gather*}
Since we know the coefficients for the similar terms in our basis, the constants $A$, $B$, $C$ can be determined from the (over-constrained) system, yielding
\[
\Tr \big(X^4\big) = \frac{3}{8} \Tr \big(M^2\big)^2 - \frac{1}{2} \Tr \big(M^4\big) - 12 \Pf (M) \] as in Theorem \ref{prop:Trafo}. In the same way, one can find the coefficients for $\Tr \big(X^6\big)$ and $\Pf (X)$, though the computations are even more lenghthy: For the Pfaffian, one first needs to find a closed formula for taking the square-root of the determinant, while for $\Tr \big(X^6\big)$ the linear system as well as the individual expressions simply increase in size.
\end{proof}

To understand the action of the order three automorphism on the basis of homogeneous invariant polynomials, note that the transformations from Theorem~\ref{prop:Trafo} are most conveniently encoded by the following matrix
 \begin{equation*}
 T=\begin{pmatrix}
1 & 0 & 0 & 0 \\
3/8 & -1/2 & -12 & 0 \\
-1/64 & 1/16 & -1/2 & 0 \\
15/64 & -15/16 & -15/2 & 1
\end{pmatrix}.
\end{equation*}
The action on a basis for homogeneous invariant polynomials of degree 6 can thus be seen as follows
\begin{equation*}
 \begin{pmatrix}
\Tr \big(\sigma(M)^2\big)^3\\ \Tr\big(\sigma(M)^2\big)\Tr\big(\sigma(M)^4\big) \\ \Tr\big(\sigma(M)^2\big) \Pf(\sigma(M)) \\ \Tr\big(\sigma(M)^6\big)
\end{pmatrix} =
\begin{pmatrix}
1 & 0 & 0 & 0 \\
3/8 & -1/2 & -12 & 0 \\
-1/64 & 1/16 & -1/2 & 0 \\
15/64 & -15/16 & -15/2 & 1
\end{pmatrix}
\begin{pmatrix}
\Tr \big(M^2\big)^3\\ \Tr\big(M^2\big)\Tr\big(M^4\big) \\ \Tr\big(M^2\big) \Pf(M) \\ \Tr\big(M^6\big)
\end{pmatrix}.
\end{equation*}
 Moreover, a reiterated action can be computed by powers of this transformation matrix $T$, through which we have
 \begin{equation*}
 T^2 = \begin{pmatrix}
1 & 0 & 0 & 0 \\
3/8 & -1/2 & 12 & 0 \\
1/64 & -1/16 & -1/2 & 0 \\
15/64 & -15/16 & 15/2 & 1
\end{pmatrix},
\end{equation*}
and $T^3 = I$ since $\sigma^3$ acts as the identity.
 Since $\sigma$ acts linearly, the fixed locus is determined by finding the eigenspace of an associated matrix, and we find the following:
 \begin{Proposition}\label{deg6}
 The space of invariant polynomials of ${\rm SO}(8,\mathbb{C})$ of degree six which are invariant under the induced action of the automorphism $\sigma$ is two-dimensional and spanned by \[\Tr\big(M^2\big)^3 \qquad \text{and}\qquad 5 \Tr\big(M^2\big) \Tr\big(M^4\big) - 8 \Tr \big(M^6\big).\]
 \end{Proposition}
 \begin{proof}
 This is verified by computing the $+1$-eigenspace of $T^t$.
 \end{proof}

The previous proposition is important because the algebra of invariant polynomials of $\mathrm{G}_2$ is generated by two homogeneous polynomials, one of degree two and one of degree six. The image of $\mathfrak{g}_2$ inside of $\mathfrak{so}(8)$ is contained in the set of matrices $M$ with eigenvalues $(0,0,\pm \eta_1, \pm \eta_2, \pm \eta_3)$ such that $\eta_1 + \eta_2 + \eta_3 = 0$, see~\cite{N3}. In terms of this representation, the two generating invariant polynomials take values
 \begin{eqnarray} c_1=\eta_1^2+\eta_2^2+\eta_3^2\qquad \text{and} \qquad c_3=(\eta_1 \eta_2 \eta_3)^2. \label{13}
 \end{eqnarray}
The following Proposition explains their role with respect to the generating set $\{p_1, p_2, p_3, p_4 \}$:

\begin{Proposition}\label{image}
The invariant polynomials $p_1$, $p_2$, $p_3$, $p_4$ of $\mathfrak{so}(8)$ restrict to invariant polynomials of~$\mathfrak{g}_2$. As invariant polynomials of $\mathfrak{g}_2$ they relate to the generating polynomials~$c_1$,~$c_3$ via
\begin{gather*}
 c_1 = \frac{1}{2} \Tr\big(M^2\big) = \frac{1}{2} p_1 , \\
 c_3 = \frac{1}{16}p_1^3 - 5 p_1 p_2 + 8 p_3.
\end{gather*}
\end{Proposition}
\begin{proof} The invariant polynomials restrict by general arguments about subgroups, namely because
\begin{equation}\label{inj}\C [\mathfrak{g}_2]^{\mathrm{G}_2} \hookrightarrow \C [\mathfrak{so}(8)]^{\mathrm{G}_2} \twoheadrightarrow \C [\mathfrak{so}(8)]^{\SO(8)}.\end{equation}
The equations are readily verified using the description of $\mathfrak{g}_2$ inside of $\mathfrak{so}(8)$ from above, since one can restrict them to matrices with eigenvalues $(0,0,\pm \eta_1, \pm \eta_2, \pm \eta_3)$ such that $\eta_1 + \eta_2 + \eta_3 = 0$.
\end{proof}

\begin{Remark}
Notice in particular that Propositions~\ref{deg6} and~\ref{image} verify that
\[
\C [\mathfrak{g}_2]^{\mathrm{G}_2} \subset \big(\C [\mathfrak{so}(8)]^{\SO (8)}\big)^\sigma
\] as expected. Moreover, one can see that for a matrix $M \in \mathfrak{g}_2 \subset \mathfrak{so}(8)$ one has $\Pf (M) = 0$ and $\Tr \big(M^4\big) = 1/2 \Tr \big(M^2\big)^2$, which shows the opposite inclusion. Although not explicitly done, this can also be deduced from \cite[Section~8.8, p.~144]{slod}. Moreover, we are very thankful to one of our reviewers who pointed out that
equality of the two invariant rings is explicitly stated and proved in \cite[Corollary~2.2.3(ii)]{review1}.
\end{Remark}

\subsection{Final remarks on further directions}\label{sec:final}
We shall conclude this short note mentioning two directions in which the present results could be useful for. However, to maintain our focus on the Lie theoretic aspect of the research, we shall leave these questions to future work.

 A natural question arising from Proposition~\ref{image} is to identify the image of $\mathcal{A}_{\mathrm{G}_2} \rightarrow \mathcal{A}_{\SO(8)}$ appearing through equation~\eqref{inj}.
When considering this question one should note that the action of $\sigma$ on the group $G=\SO(8)$ requires a choice of splitting of the sequence $0 \rightarrow \operatorname{Inn}(G) \rightarrow \operatorname{Aut}(G) \rightarrow \operatorname{Out}(G) \rightarrow 0$.
 This sequence is always split but not canonically so: A choice of splitting is equivalent to a choice of Cartan and Borel for $G$. Further sources to investigate this direction appear in \cite{anton, oscar}, and references therein. The action on $\mathcal{M}_{G_\mathbb{C}}$ is independent of choices, since any two representatives differ by conjugation, which via non-abelian Hodge theory acts trivially on $\mathcal{M}_{G_\mathbb{C}}$.

 Finally, with views towards applications within Langlands duality and mirror symmetry, it is also natural to ask what the effect of triality is on Lagrangian subspaces of the moduli space of Higgs bundles defined through other automorphisms, such as those used in \cite{slices,cmc}.
In this direction, the reader might find of interest the work in \cite[Section~10]{oscar} where the authors show how the triality automorphism moves in a cyclic way the moduli spaces of Higgs bundles corresponding to three different realizations of the real forms ${\rm SO}_0(3,5)$ and ${\rm SO}_0(1,7)$.

\subsection*{Acknowledgements}
The authors are thankful to S.~Rayan for his thorough comments on a draft of the manuscript. The work of S.S.\ is partially supported by NSF grants DMS 1107452, 1107263, 1107367 ``RNMS: GEometric structures And Representation varieties (the GEAR Network)''.
 L.P.S.\ was partially supported by NSF DMS 1509693 and NSF CAREER Award DMS 1749013, as well as by the Alexander Von Humboldt foundation. This material is also based upon work supported by NSF DMS 1440140 while L.P.S.\ was in residence at the Mathematical Sciences
Research Institute in Berkeley, California, during the Fall 2019 semester. Both authors are thankful for the support of the Simons Center for Geometry and Physics during the Spring 2019 program on {\it Geometry and Physics of Hitchin systems}.

\pdfbookmark[1]{References}{ref}
\LastPageEnding

\end{document}